\def\F{\mathbb F}

  \def\a{\alpha}
  \def\b{\beta}
  \def\g{\gamma}

  \def\D{\Delta}
  \def\G{\Gamma}
  \def\e{\epsilon}

  \def\la{\langle}
  \def\ra{\rangle}
  \def\s{\sigma}
  \def\o{\omega}
  \def\O{\Omega} 
 
  \def\pf{{\it Proof. }$\;\;$}
  \def\no{\noindent}
  \def\hal{\unskip\nobreak\hfil\penalty50\hskip10pt\hbox{}\nobreak
  \hfill\vrule height 5pt width 6pt depth 1pt\par\vskip 2mm}

    \documentclass[11pt,a4paper]{article}
            \usepackage{amsfonts,amssymb}
\usepackage{amsmath}

\begin{document}

  \title{The classification of $\frac{3}{2}$-transitive permutation groups and $\frac{1}{2}$-transitive linear groups}
  \author{ Martin W. Liebeck, Cheryl E. Praeger and Jan Saxl}
  \maketitle

 \begin{abstract}
A linear group $G\le GL(V)$, where $V$ is a finite vector space, is called $\frac{1}{2}$-transitive if all the $G$-orbits on the set of nonzero vectors have the same size. We complete the classification of all the $\frac{1}{2}$-transitive linear groups. As a consequence we complete the determination of the finite $\frac{3}{2}$-transitive permutation groups -- the transitive groups for which a point-stabilizer has all its nontrivial orbits of the same size. We also determine the $(k+\frac{1}{2})$-transitive groups for integers $k\ge 2$.
 \end{abstract}


\newtheorem{theorem}{Theorem}
  \newtheorem{thm}{Theorem}[section]
  \newtheorem{prop}[thm]{Proposition}
  \newtheorem{lem}[thm]{Lemma}
  \newtheorem{cor}[theorem]{Corollary}
\newtheorem{proposition}[theorem]{Proposition}

\section{Introduction}\label{one}
The concept of a finite {\it $\frac{3}{2}$-transitive} permutation group -- a non-regular transitive group in which all the nontrivial orbits of a point-stabilizer have equal size -- was introduced by Wielandt in his book \cite[\S 10]{W}. Examples are 2-transitive groups and Frobenius groups: for the former, a point-stabilizer has just one nontrivial orbit, and for the latter, every nontrivial orbit of a point-stabilizer is regular. Further examples are provided by normal subgroups of 2-transitive groups; indeed, one of the reasons for Wielandt's definition was that normal subgroups of 2-transitive groups are necessarily $\frac{3}{2}$-transitive.

Wielandt proved that any $\frac{3}{2}$-transitive group is either primitive or a Frobenius group (\cite[Theorem 10.4]{W}). Following this, a substantial study of $\frac{3}{2}$-transitive groups was undertaken by Passman in \cite{pass1,pass2}, in particular completely determining the soluble examples. More recent steps towards the classification of the primitive $\frac{3}{2}$-transitive groups were taken in \cite{3half} and \cite{arith}. In \cite{3half} it was proved that 
primitive $\frac{3}{2}$-transitive groups are either affine or almost simple, and the almost simple examples were determined. For the affine case, consider an affine group $T(V)G \le AGL(V)$, where $V$ is a finite vector space, $T(V)$ is the group of translations, and $G \le GL(V)$; this group is $\frac{3}{2}$-transitive if and only if the linear group $G$ is {\it $\frac{1}{2}$-transitive} -- that is, all the orbits of $G$ on the set $V^\sharp$ of nonzero vectors have the same size. 
The $\frac{1}{2}$-transitive linear groups of order divisible by $p$ (the characteristic of the field over which $V$ is defined) were determined in \cite[Theorem 6]{arith}. 

The main result of this paper completes the classification of $\frac{1}{2}$-transitive linear groups. In the statement, by a {\it semiregular} group, we mean a permutation group all of whose orbits are regular.

\begin{theorem}\label{hope}
Let $G \le GL(V) = GL_d(p)$ ($p$ prime) be an insoluble $p'$-group, and suppose $G$ is $\frac{1}{2}$-transitive on $V^\sharp$. Then one of the following holds:
\begin{itemize}
\item[(i)] $G$ is semiregular on $V^\sharp$;
\item[(ii)] $d=2$, $p=11,19$ or $29$, and $SL_2(5)\triangleleft G \le GL_2(p)$;
\item[(iii)] $d=4$, $p=13$, and $SL_2(5) \triangleleft G \le \G L_2(p^2) \le GL_4(p)$.
\end{itemize}
In (ii) and (iii), the non-semiregular possibilities for $G$ are given in Table $\ref{possg}$.
\end{theorem}

\begin{table}\label{possg}
\caption{Orbit sizes of $\frac{1}{2}$-transitive groups in Theorem \ref{hope}(ii),(iii)}
\[
\begin{array}{|l|l|l|l|}
\hline
p^d & |G| & \hbox{orbit size on }V^\sharp & \hbox{number of orbits} \\
\hline
11^2 & 600 & 120 & 1 \\
19^2 & 360 & 120 & 3 \\
         &  1080 & 360 & 1 \\
29^2 & 240 & 120 & 7 \\
         & 1680 & 840 & 1 \\
13^4 & 3360 & 1680 & 17 \\
\hline
\end{array}
\]
\end{table}

\no {\bf Remarks } 1. In conclusion (i) of the theorem, the corresponding affine permutation group $T(V)G$ (acting on $V$) is a Frobenius group, and $G$ is a Frobenius complement (see Proposition \ref{frob} for the structure of these).

\vspace{2mm}
\no 2. In conclusion (ii), $\F_p^*R$ acts transitively on $V^\sharp$, where $R = SL_2(5)$ and $\F_p^*$ is the group of scalars in $GL(V)$, and $G = Z_0R$ for some $Z_0 \le \F_p^*$. Here $G \triangleleft \F_p^*R$ (hence is  $\frac{1}{2}$-transitive, since in general, a normal subgroup of a transitive group is $\frac{1}{2}$-transitive).

\no 3. The $\frac{1}{2}$-transitive group $G$ in part (iii) is more interesting. Here $G = (Z_0R).2 \le \G L_2(13^2)$, where
$R = SL_2(5)$ and $Z_0$ is a subgroup of $\F_{13^2}^*$ of order 28, and $G\cap GL_2(13^2) = Z_0R$ has orbits on 1-spaces of sizes 20, 30, 60, 60.

\vspace{4mm}
Combining Theorem \ref{hope} with the soluble case in \cite{pass1,pass2} and the $p$-modular case in \cite[Theorem 6]{arith}, we have the following classification of $\frac{1}{2}$-transitive linear groups. In the statement, for $q$ an odd prime power, $S_0(q)$ is the subgroup of $GL_2(q)$ of order $4(q-1)$ consisting of all monomial matrices of determinant $\pm 1$.

\begin{cor}\label{allthree}
If $G \le GL(V) = GL_d(p)$ is $\frac{1}{2}$-transitive on $V^\sharp$, then one of the following holds:
\begin{itemize}
\item[(i)] $G$ is transitive on $V^\sharp$;
\item[(ii)] $G \le \G L_1(p^d)$;
\item[(iii)] $G$ is a Frobenius complement acting semiregularly on $V^\sharp$;
\item[(iv)] $G = S_0(p^{d/2})$ with $p$ odd;
\item[(v)] $G$ is soluble and $p^d = 3^2,5^2,7^2,11^2,17^2$ or $3^4$;
\item[(vi)] $SL_2(5) \triangleleft G \le \G L_2(p^{d/2})$, where $p^{d/2} = 9$, $11$, $19$, $29$ or $169$.
\end{itemize}
\end{cor}

The classification of $\frac{3}{2}$-transitive permutation groups follows immediately from this result and those in \cite{3half}. For completeness, we state it here.

\begin{cor}\label{32classn}
Let $X$ be a $\frac{3}{2}$-transitive permutation group of degree $n$. Then one of the following holds:
\begin{itemize}
\item[(i)] $X$ is $2$-transitive;
\item[(ii)] $X$ is a Frobenius group;
\item[(iii)] $X$ is affine: $X = T(V)G \le AGL(V)$, where $G \le GL(V)$ is a $\frac{1}{2}$-transitive linear group, given by Corollary $\ref{allthree}$;
\item[(iv)] $X$ is almost simple: either
\begin{itemize}
\item[(a)] $n=21$, $X = A_7$ or $S_7$ acting on the set of pairs in $\{1,\ldots ,7\}$, or
\item[(b)] $n=\frac{1}{2}q(q-1)$ where $q=2^f\ge 8$, and either $G = PSL_2(q)$, or $G = P\G L_2(q)$ with $f$ prime.
\end{itemize}
\end{itemize}
\end{cor}

Turning to higher transitivity, recall (again from \cite{W}) that for a positive integer $k$, a permutation group is $(k+\frac{1}{2})$-transitive if it is $k$-transitive and the stabilizer of $k$ points has orbits of equal size on the remaining points. For $k\ge 2$ such groups are of course 2-transitive so belong to the known  list of such groups. Nevertheless, their classification has some interesting features and we record this in the following result.

\begin{proposition}\label{khalf}
Let $k\ge 2$ be an integer, and let $X$ be a  $(k+\frac{1}{2})$-transitive permutation group of degree $n\ge k+1$. Then one of the following holds:
\begin{itemize}
\item[(i)] $X$ is $(k+1)$-transitive;
\item[(ii)] $X$ is sharply $k$-transitive;
\item[(iii)] $k=3$ and $X = P\Gamma L_2(2^p)$ with $p$ an odd prime, of degree $2^p+1$;
\item[(iv)] $k=2$ and one of:
\begin{itemize}
\item[] $L_2(q) \triangleleft X \le P \Gamma  L_2(q)$ of degree $q+1$;
\item[] $X = Sz(q)$, a Suzuki group of degree $q^2+1$;
\item[] $X = A\Gamma L_1(2^p)$ with $p$ prime, of degree $2^p$.
\end{itemize}
\end{itemize}
\end{proposition}

\no {\bf Remarks } 1. The sharply $k$-transitive groups were classified by Jordan for $k \ge 4$ and by Zassenhaus for $k=2$ or 3; see \cite[\S 7.6]{DM}.

\vspace{2mm}
\no 2. In conclusion (iv), the groups $Sz(q)$ and $A\Gamma L_1(2^p)$ are Zassenhaus groups -- that is, 2-transitive groups in which all 3-point stabilizers are trivial (so that all orbits of a 2-point stabilizer are regular). The groups $X$ with socle $L_2(q)$ are all $\frac{5}{2}$-transitive, being normal subgroups of the 3-transitive group $P\Gamma L_2(q)$; some are 3-transitive, some are Zassenhaus groups, and some are neither.

\vspace{4mm} 
The paper consists of two further sections, one proving Theorem \ref{hope}, and the other Proposition \ref{khalf}.

\vspace{4mm}
\no{\bf  Acknowledgements } We thank Eamonn O'Brien for assistance with the Magma computations in the paper. 
The second author acknowledges the support of Australian Research Council Discovery Project Grant DP140100416.

\section{Proof of Theorem \ref{hope}}

Throughout the proof, we shall use the following well-known result about the structure of Frobenius complements, due to Zassenhaus.

\begin{prop}\label{frob} {\rm (\cite[Theorem 18.6]{passbk})} Let $G$ be a Frobenius complement.
\begin{itemize}
\item[(i)] The Sylow subgroups of $G$ are cyclic or generalized quaternion.
\item[(ii)] If $G$ is insoluble, then it has a subgroup of index 1 or 2 of the form $SL_2(5) \times Z$, where $Z$ is a group of order coprime to $30$, all of whose Sylow subgroups are cyclic.
\end{itemize}
\end{prop}

The following result is important in our inductive proof of Theorem \ref{hope}.

\begin{prop}\label{sl25}
Let $R = SL_2(5)$, let $p>5$ be a prime, and let $V$ be a nontrivial absolutely irreducible $\F_qR$-module, where $q=p^a$. Regard $R$ as a subgroup of $GL(V)$, and let $G$ be a group such that $R\triangleleft G \le \G L(V)$.
\begin{itemize}
\item[(i)] If $R$ is semiregular on $V^\sharp$, then $\dim V = 2$.
\item[(ii)] Suppose $\dim V = 2$ and $G$ has no regular orbit on the set $P_1(V)$ of
 $1$-spaces in $V$. Then either $q \in \{p,p^2\}$ with $p\le 61$, or $q = 7^4$.
\item[(iii)] If $\dim V=2$ and $G$ is $\frac{1}{2}$-transitive but not semiregular on $V^\sharp$, then $q=11,19,29$ or $169$. Conversely, for each of these values of $q$ there are examples of $\frac{1}{2}$-transitive, non-semiregular groups $G$, and they are as in Table $\ref{possg}$ of Theorem $\ref{hope}$.
\end{itemize}
\end{prop}

\pf (i) The irreducible $R$-modules and their Brauer characters can be found in \cite{Atlas}, and have dimensions 2, 3, 4, 5 or 6. For those of dimension 3 or 5, the acting group is $R/Z(R) \cong A_5$, and involutions fix nonzero vectors; and for those of dimension 4 or 6, elements of order 3 fix vectors.

(ii) Let $\dim V = 2$, and suppose $G$ has no regular orbit on $P_1(V)$. Assume for a contradiction that $q$ is not as in the  conclusion of (ii). In particular, $q>61$ (recall that $p>5$).

Write $\bar R = R/Z(R) \cong A_5$ and $\bar G = G/(G\cap \F_q^*)$. Now $N_{PGL(V)}(\bar R) = \bar R$, so it follows that $\bar G = \bar R\la \s\ra$ for some $\s \in P\G L(V)$ (possibly trivial). Note that if $p\equiv \pm 2\hbox{ mod }5$ then $\F_{p^2} \subseteq \F_q$.

Consider the action of $\bar R\cong A_5$ on $P_1(V)$. As $A_5$ has 31 nontrivial cyclic subgroups, and each of these fixes at most two 1-spaces, it follows that $\bar R$ has at least $(q-62)/60$ regular orbits on $P_1(V)$. Since $q>61$, $\bar R$ has a regular orbit, and so $\bar G\ne \bar R$ by our assumption.

Let $r$ be the order of the element $\s$ modulo $\bar R$ (so that $\F_{p^r}\subseteq \F_q$). If there is a regular $\bar R$-orbit $\D_0$ on $P_1(V )$ that is not fixed by $\s^i$ for any $i$ with $1\le i\le r-1$, then $\bar G_{\D_0} = \bar R$ and so $\bar G_{\la v \ra} = 1$ for $\la v\ra \in \D_0$ and $G$ has a regular orbit on $P_1(V)$, a contradiction. Hence $r>1$, and for each regular $\bar R$-orbit $\D$, there is a subgroup 
$\la \s^{i(\D)}\ra$, of prime order modulo $\bar R$, which fixes $\D$ setwise. Moreover, for $\la v\ra \in \D$, there exists $x \in \bar R$ such that $x \s^{i(\D)}$ fixes $\la v\ra$. Since there are at least $q-62$ elements of $P_1(V)$ in regular $\bar R$-orbits, 
it follows that 
\begin{equation}\label{ineq}
|\bigcup {\rm fix}_{P_1(V)}(x\s^{j})| \ge q-62,
\end{equation}
where the union is over all $x \in \bar R$ and all $j$ dividing $r$ with $r/j$ prime. Let $s = r/j$ for such $j$, and let $x \in \bar R$. If $(x\s^j)^s \ne 1$ then $(x\s^j)^s \in \bar R$ fixes at most two 1-spaces, and so 
$|{\rm fix}(x\s^j)| \le 2$; and if $(x\s^j)^s = 1$, then $x\s^j$ is $PGL(V)$-conjugate to a field automorphism of order $s$, and $|{\rm fix}(x\s^j)| = q^{1/s}+1$. Hence (\ref{ineq}) implies that 
\begin{equation}\label{ineq2}
60\sum_{s|r, s\;prime}(q^{1/s}+1) \ge q-62.
\end{equation}
Recall that $p>5$ and $\F_{p^r}\subseteq \F_q$.

Suppose that $6|r$. The terms in the sum on the left hand side of (\ref{ineq2}) with $s\ge 5$ add to at most $r(q^{1/5}+1)$, which is easily seen to be less than $q^{1/2}+1$. Hence (\ref{ineq2}) gives
\[
2(q^{1/2}+1)+(q^{1/3}+1) \ge \frac{q-62}{60}.
\]
Putting $y=q^{1/6}$ this yields $120y^3+60y^2+242 \ge y^6$, which is false for $y \ge 7$. Similarly, when ${\rm hcf}(r,6) = 1$ or 3, we find that (\ref{ineq2}) fails. Consequently ${\rm hcf}(r,6) = 2$, and (\ref{ineq2}) gives $2(q^{1/2}+1) \ge (q-62)/60$, which implies that $q^{1/2}\le 121$. Hence (as $p>5$ and $q=p^a$ with $a$ even), either $q=p^2$ or $q=7^4$ or $11^4$. Then further use of 
(\ref{ineq2}) gives $p\le 61$ in the former case, and also shows that $q\ne 11^4$.  But now we have shown that $q$ is as in (ii), contrary to assumption. This completes the proof. 

(iii) Suppose $G$ is $\frac{1}{2}$-transitive but not semiregular on $V^\sharp$.
If $G$ has a regular orbit on $P_1(V)$, then it has a regular orbit on $V^\sharp$, which is not possible by the assumption in the previous sentence. Hence  $q$ must be as in the conclusion of part (ii). For these values of $q$, we use Magma \cite{Mag} to construct $R \cong SL_2(5)$ in $SL_2(q)$, and for all subgroups of $\G L_2(q)$ normalizing $R$, compute whether they are $\frac{1}{2}$-transitive and non-semiregular. We find that such groups exist precisely when $q$ is 11, 19, 29 or 169, and the examples are as in Table \ref{possg}. \hal

\vspace{4mm}
Note that part (ii) of the proposition follows from \cite[Theorem 2.2]{KP} in the case where $R$ is $\F_p$-irreducible on $V$. We shall need the more general case proved above.

\vspace{2mm}
We now embark on the proof of Theorem \ref{hope}. Suppose that $G$ is a minimal counterexample. That is, 
\begin{itemize}
\item $G \le GL_d(p) = GL(V)$ is an insoluble, $\frac{1}{2}$-transitive $p'$-group, 
\item $G$ is not semiregular on $V^\sharp$, and $G$ is not as in (ii) or (iii) of the theorem, and 
\item $G$ is minimal subject to these conditions. 
\end{itemize}
Observe that since $G$ is  $\frac{1}{2}$-transitive and not semiregular, it cannot have a regular orbit on $V$.

The affine permutation group $VG \le AGL(V)$ is $\frac{3}{2}$-transitive on $V$ and not a Frobenius group, hence is primitive by \cite[Theorem 10.4]{W}. It follows that $G$ is irreducible on $V$.

By \cite[Theorem 1.1]{pass2}, $G$ acts primitively as a linear group on $V$. Choose $q=p^k$ maximal such that $G \le \G L_n(q) \le GL_d(p)$, where $d = nk$. Write $V = V_n(q)$, $G_0 = G \cap GL_n(q)$, $K = \F_q$ and $Z = G_0 \cap K^*$, the group of scalars in $G_0$. Since $G$ is insoluble, $n \ge 2$. Also $G_0$ is absolutely irreducible on $V$ (see \cite[Lemma 12.1]{arith}), so $Z = Z(G_0)$.

\begin{lem}\label{norm1}
Let $N$ be a normal subgroup of $G$ with $N \le G_0$ and $N \not \le Z$, and let $U$ be an irreducible $KN$-submodule of $V$. Then the following hold:
\begin{itemize}
\item[(i)] $N$ acts faithfully and absolutely irreducibly on $U$;
\item[(ii)] $N$ is not cyclic;
\item[(iii)] $G_U$ acts $\frac{1}{2}$-transitively on $U^\sharp$;
\item[(iv)] if $(G_U)^U$ is insoluble and not semiregular, and $(N^{(\infty)},|U|) \ne (SL_2(5),q^2)$ with $q\in \{11,19,29,169\}$, then $U=V$.
\end{itemize}
\end{lem}

\pf As $G$ is primitive on $V$, Clifford's theorem implies that $V\downarrow N$ is homogeneous, so that $V\downarrow N = 
U \oplus U_2 \oplus \cdots \oplus U_r$ with each $U_i \cong U$. Hence $N$ is faithful on $U$; it is also absolutely irreducible, as in the proof of \cite[Lemma 12.2]{arith}. Hence (i) holds, and (ii) follows.

 To see (iii), let $v \in U^\sharp$, $n \in N$ and $g \in G_v$. Then $vng = vgn' = vn'$ for some $n' \in N$. Hence $\{vn : n\in N\}$ is invariant under $G_v$. As $U$ is irreducible under $N$, $\{vn:n\in N\}$ spans $U$, and hence $G_v$ stabilises $U$. Therefore 
\[
|G:G_v| = |G:G_U|\cdot |G_U:G_v|.
\]
As $G$ is $\frac{1}{2}$-transitive this is independent of $v \in U^\sharp$, and hence $G_U$ is $\frac{1}{2}$-transitive on 
$U^\sharp$, as in (iii).  

Finally, (iv) follows by the minimality of $G$. \hal

\vspace{4mm}
By \cite[Theorem A]{pass2}, $O_r(G_0)$ is cyclic for each odd prime $r$, and hence is central by Lemma \ref{norm1}(ii). 
Consequently $F(G_0) = ZE$ where $E = O_2(G_0)$. Moreover \cite[Theorem A]{pass2} also shows that $\Phi(E)$
is cyclic, hence contained in $Z$, and $|E/\Phi(E)| \le 2^8$.

Now let $F^*(G_0) = ZER_1\cdots R_k$, a commuting product with each $R_i$ quasisimple (possibly $k=0$).

\begin{lem}\label{k0} We have $k\ge 1$.
\end{lem}

\pf Suppose $k=0$, and write $N = F^*(G_0) = ZE$. 
Since $V\downarrow G$ is primitive, every characteristic abelian subgroup of $E$ is cyclic, so $E$ is a 2-group of symplectic type. By a result of Philip Hall (\cite[23.9]{aschbk}), we have $E = E_1 \circ F$ where $E_1$ is either 1 or extraspecial, and $F$ is cyclic, dihedral, semidihedral or generalised quaternion; in the latter three cases, $|F|\ge 16$.  
Since $N = F^*(G_0)$ we have $C_{G_0}(N) \le N$ and $G_0/C_{G_0}(N) \le {\rm Aut}(N)$. Hence ${\rm Aut}(N)$ must be insoluble, and it follows that $|E_1/\Phi(E_1)|\ge 2^4$. 

Now $E$ has a characteristic subgroup $E_0 = E_1\circ L$, where $L = C_4$ if 4 divides $|F|$ and $L=1$ otherwise. Then $E_0 \triangleleft G$. 
 Let $U$ be an irreducible $KE_0$-submodule of $V$. By Lemma \ref{norm1}, $E_0$ is faithful on $U$ and $G_U$ is $\frac{1}{2}$-transitive on $U^\sharp$. Write $H =  (G_U)^U$.

Assume that $H$ is soluble. As $H$ is $\frac{1}{2}$-transitive on $U^\sharp$, it is therefore given by \cite[Theorem B]{pass2}, which implies that one of the following holds:
\begin{itemize}
\item[(a)] $H$ is a Frobenius complement;
\item[(b)] $H \le \G L_1(q^u)$, where $|U| = q^u$;
\item[(c)] $H \le GL_2(q^u)$ with $|U| = q^{2u}$, and $H$ consists of all monomial matrices of determinant $\pm 1$;
\item[(d)] $|U| = p^2$ with $p \in \{3,5,7,11,17\}$, or $|U| = 3^4$.
\end{itemize}
In all cases except the last one in (d), it follows (using Proposition \ref{frob}(i) for (a)) that $|E_0/\Phi(E_0)|\le 2^2$, which is a contradiction. In the exceptional case $|U| = 3^4$ and $|E_0/\Phi(_0E)|= 2^4$.  
But in this case any $3'$-subgroup of ${\rm Aut}(N)$ is soluble, and hence $G_0$ is soluble, again a contradiction. 

Hence $H$ is insoluble.  As $H$ is not a Frobenius complement by Proposition \ref{frob}(ii), it is not semiregular on $U^\sharp$, 
and so Lemma \ref{norm1}(iv) implies that $U=V$. Hence $E_0$ is irreducible on $V$ and so $F$ is cyclic and $N = ZE = ZE_0$. We have $|E_0/\Phi(E_0)| \le 2^8$ by \cite[Theorem A]{pass2}, and hence $|E_0/\Phi(E_0)| = 2^{2m}$ with $m=2,3$ or 4.

 \vspace{2mm} 
\no {\bf Case $m=4$.}
Suppose first that $m=4$, so $E_1 = 2^{1+8}$ and $\dim V = 16$. By \cite[Lemmas 2.6, 2.10]{pass2} we have $E_1=E_0$, so that $|Z|_2 = 2$ and 
$G_0 \le Z \circ 2^{1+8}.O_8^\e (2)$ ($\e = \pm$). Also \cite[Lemma 2.4]{pass2} gives $(p^2-1)_2 \ge 2^4$, hence $p\ge 7$, and the proof of \cite[Lemma 2.12]{pass2} gives $|G/N| \ge q^8/2^9$. Since $G/N \le O_8^\e (2)$, it follows that $q=7$. Hence $G/N$ is an insoluble 
$7'$-subgroup of $O_8^\e (2)$ of order greater than $7^8/2^9$. Using \cite{Atlas}, we see that such a subgroup is contained in one of the following subgroups of $O_8^\e(2)$:
\[
\begin{array}{l}
2^6.O_6^-(2),\,2^{1+8}.(S_3\times S_5)\; (\e=-) \\
S_3\times O_6^-(2),\,2^6.(S_6\times 2),\,2^6.(S_5\times S_3),\,(S_5\times S_5).2 \; (\e=+)
\end{array}
\]

We now consider elements of order 3 in $G$. These are elements $t_k$ lying in subgroups $O_2^-(2)^k$ of $O_8^\e (2)$ for $1\le k\le 4$ and acting on the 16-dimensional space $V$ as a tensor product of $k$ diagonal matrices $(\o, \o^{-1})$ with an identity matrix of dimension $2^{4-k}$, where $\o \in K^*$ is a primitive cube root of 1; there are also scalar multiples $\o t_k$ if $Z$ contains $\o I$. We compute the action of $t_k$ on $V$ and also the class of the image of $t_k$ in $O_8^\e (2)$ in Atlas notation, as follows:
\[
\begin{array}{l|l|l}
k & \hbox{action of }t_k \hbox{ on }V & \hbox{Atlas notation} \\
\hline
1 & (\o^{(8)},\o^{-1\,(8)}) & 3A\,(\e=-),\,3A\,(\e=+) \\
2 & (1^{(8)},\o^{(4)},\o^{-1\,(4)}) & 3B\,(\e=-),\,3E\,(\e=+) \\
3 &  (1^{(4)},\o^{(6)},\o^{-1\,(6)}) & 3C\,(\e=-),\,3D\,(\e=+) \\
4 &  (1^{(6)},\o^{(5)},\o^{-1\,(5)}) & -\,(\e=-),\,3BC\,(\e=+) \\
\hline
\end{array}
\]
Hence every element of order 3 in $G$ has fixed point space on $V$ of dimension at most 8. Considering the above subgroups of $O_8^\e(2)$, we compute that the total number of elements of order 3 in $G$ is less than $2^{20}$. If $G$ contains an element of order 3 fixing a nonzero vector in $V$, then as $G$ is $\frac{1}{2}$-transitive, every nonzero vector is fixed by some element of $G$ of order 3. Hence $V$ is the union of the subspaces $C_V(t)$ over $t\in G$ of order 3, so that 
\begin{equation}\label{3elt}
|V| \le \sum_{t\in G, |t|=3}|C_V(t)|.
\end{equation}
This yields $7^{16} < 2^{20}\cdot 7^8$, which is false. 

It follows that $G$ contains no element of order 3 fixing a nonzero vector. So every element of order 3 in $G/N$ is conjugate to $t_1$. 

We now complete the argument by considering involutions in $G$. Now $G$ certainly contains involutions which fix nonzero vectors, so arguing as above we have
\begin{equation}\label{invol}
|V| \le \sum_{t\in G, |t|=2}|C_V(t)|.
\end{equation}
The group $G/N$ is an insoluble $7'$-subgroup of $O_8^\e(2)$, all of whose elements of order 3 are conjugates of $t_1$. Using Magma \cite{Mag}, we compute that there are 206 such subgroups if $\e=+$, and 59 if $\e=-$. For each of these possibilities for $G/N$ we compute the list of involutions of $G$ and their fixed point space dimensions. All possibilities contradict (\ref{invol}). For example, when $\e = -$ the largest possibility for $G$ has 188 involutions with fixed space of dimension 12; 74886 with dimension 8; and 188 with dimension 4. Hence (\ref{invol}) gives
\[
7^{16} \le 188\cdot (7^{12}+7^4) + 74886\cdot 7^8,
\]
which is false. This completes the proof for $m=4$.

 \vspace{2mm} 
\no {\bf Case $m=3$. }
Now suppose $m=3$, so that $\dim V=8$. This case is handled along similar lines to the previous one. By \cite[Lemma 2.9]{pass2}, either $|Z|_2 = 2$ and $G_0/N \le O^\e_6(2)$, or 4 divides $|Z|$ and $G$ contains a field automorphism of order 2 (so that $q$ is a square), and $G_0/N \le Sp_6(2)$. As $G_0$ is insoluble, its order is divisible by 2 and 3, so $p\ge 5$. Also each non-central 
involution in $G_0$ fixes a nonzero vector.

Assume now that 7 divides $|G|$. If 7 divides $|G/G_0|$ then $q \ge 5^7$ and we easily obtain a contradiction using (\ref{invol}); so 7 divides $|G_0|$. Elements of order 7 in $G_0$ act on $V$ as $(1^2,\o,\o^2,\ldots ,\o^6)$ where $\o$ is a 7th root of 1 in the algebraic closure of $\F_q$ (since they are rational in $O_6^+(2)$). In particular they fix nonzero vectors, so $\frac{1}{2}$-transitivity implies
\begin{equation}\label{7elt}
|V| \le \sum_{t\in G, |t|=7}|C_V(t)|.
\end{equation}
The number of elements of order 7 in $Sp_6(2)$ is 207360, and hence the number in $G_0$ is at most 
$(q-1,7)\cdot 2^6\cdot 207360$. Each fixes at most $q^2$ vectors, so (\ref{7elt}) gives
\[
q^8 \le (q-1,7)\cdot 2^6\cdot 207360\cdot q^2,
\]
which implies that $q\le 13$. Hence $q = 5,11$ or 13 (not 7 as $G_0$ is a $p'$-group). As $q$ is prime, by the first observation in this case, we have $|Z|_2 = 2$ and $G/N \le O^+_6(2)$. But then the number of elements of order 7 in $G$ is at most 
$2^6\cdot 5760$, so (\ref{7elt}) forces $q=5$. So $G/N$ is an insoluble $5'$-subgroup of $O_6^+(2)$, and now we use Magma to see that such a group $G$ is not $\frac{1}{2}$-transitive on the nonzero vectors of $V = V_8(5)$.

Therefore 7 does not divide $|G|$. It follows that $G_0/N$ is contained in one of the following subgroups of $Sp_6(2)$:
\[
O_6^-(2),\,S_6\times S_3,\,2^5.S_6.
\]
As $G_0$ is insoluble and a $p'$-group, we have $p\ge 7$. We now consider elements of order 3 in $G$. 
These are conjugate to elements $t_k$ ($1\le k\le 3$) lying in subgroups $(O_2^-(2))^k$ of $Sp_6(2)$, and acting on $V$ as follows: 
\[
\begin{array}{l}
t_1: \,(\o^{(4)},\o^{-1\,(4)}), \\ 
t_2: \,(1^4,\o^{(2)},\o^{-1\,(2)}), \\
t_3: \,(1^2,\o^{(3)},\o^{-1\,(3)}).
\end{array}
\]
Suppose $G$ has an element of order 3 which fixes nonzero vectors in $V$, so that (\ref{3elt}) holds. We argue as in the previous case that $q$ is not a cube, so 3 does not divide $|G/G_0|$. In $O_6^-(2)$, the numbers of elements conjugate to $t_1,t_2,t_3$ are 240, 480, 80 respectively. Hence, if $G_0/N \le O_6^-(2)$ then (\ref{3elt}) gives
\[
q^8 \le 2^4\cdot 480q^4+2^6\cdot 80q^2 +2^3\cdot 240q^4+2^5\cdot 480q^2+2^7\cdot 80q^3
\]
where the last three terms are only present if 3 divides $|Z|$. This gives $q=7$. Similarly $q=7$ is the only possibility if $G_0/N$ is 
contained in $S_6\times S_3$ or $2^5.S_6$. But now we compute using Magma that such groups $G$ are not $\frac{1}{2}$-transitive on the nonzero vectors of $V = V_8(7)$.

Thus all elements of order 3 in $G$ are fixed point free on $V^\sharp$, and hence 
$G_0/N$ is an insoluble $7'$-subgroup of $Sp_6(2)$, all of whose elements of order 3 are conjugate to $t_1$. We compute that there are 10 such subgroups, and for each of them, (\ref{invol}) implies that $q=7$ is the only possibility: for example, the largest possible $G_0$ has 60 (resp. 3526, 60) involutions with fixed point spaces on $V$ of dimension 6 (resp. 4, 2), so (\ref{invol}) yields
\[
q^8 \le 60q^6+3526q^4+60q^2,
\]
hence $q=7$. Finally, we compute that none of the possible subgroups $G$ is $\frac{1}{2}$-transitive on the nonzero vectors of $V = V_8(7)$. 

 \vspace{2mm} 
\no {\bf Case $m=2$. }
Now suppose $m=2$, so that $\dim V = 4$. Then $G_0/N$ is an insoluble subgroup of $Sp_4(2)$, so is isomorphic to $S_6,A_6,S_5$ or $A_5$.

Assume that $G_0/N$ is $A_6$ or $S_6$. Then 4 divides $|Z|$ (so divides $q-1$). Elements of $G_0$ of order 3 are conjugate to $t_k$ ($k=1,2$) lying in $Sp_2(2)^k$; and $t_1$ acts on $V$ as $(\o^{(2)},\o^{-1\,(2)})$, $t_2$ as $(1^2,\o,\o^{-1})$. 
By assumption $G_0$ contains elements in both classes, so (\ref{3elt}) yields
\[
q^4 \le 2^4\cdot 40q^2 +2\cdot 2^4\cdot 40q+2\cdot 2^2\cdot 40q^2,
\]
where the last two terms are present only if 3 divides $|Z|$ (hence also $q-1$). Since 4 divides $q-1$, we conclude that 
$q = 13$ or 17 in this case.

Now assume $G_0/N$ is $A_5$ or $S_5$. As $G$ is a $p'$-group, $p\ge 7$. We compute that $G_0$ has at most 230 involutions, so (\ref{invol}) gives 
$q^4 \le 230q^2$, whence $q\le 13$. 

Thus in all cases, we have $q = 7, 11, 13$ or 17. We now compute that none of the possibilities for $G$ is $\frac{1}{2}$-transitive on the nonzero vectors of $V = V_4(q)$. This completes the proof of the lemma. \hal

\begin{lem}\label{three}
Either $|E/\Phi(E)|\le 2^2$, or  $|E/\Phi(E)|= 2^4$ and $p=3$.
\end{lem}

\pf The result is trivial if $E\le Z$, so suppose is not the case. Let $N = ZE \triangleleft G$, and let $U$ be an irreducible $KN$-submodule of $V$. By Lemma \ref{norm1}, $N$ is faithful on $U$ and $G_U$ is $\frac{1}{2}$-transitive on $U^\sharp$. Write $H =  (G_U)^U$.

Assume first that $H$ is insoluble. Now $H$ is not semiregular on $U^\sharp$ (as it is not a Frobenius complement by Proposition \ref{frob}, having $N\cong N^U$ as a normal subgroup), so Lemma \ref{norm1}(iv) implies that $U=V$. But then $N = ZE$ is irreducible on $V$, which forces $k=0$, contrary to Lemma \ref{k0}.

Hence  $H$ is soluble. As it is $\frac{1}{2}$-transitive on $U^\sharp$, it is therefore given by \cite[Theorem B]{pass2}; the list is given under (a)-(d) in the proof of Lemma \ref{k0}. In all cases except the last one in (d), it follows that $|E/\Phi(E)|\le 2^2$; in the exceptional case $|U| = 3^4$ and 
$|E/\Phi(E)|= 2^4$. Hence the conclusion of the lemma holds. \hal

\begin{lem}\label{rinorm}
If $R_i \triangleleft G$, then $R_i = SL_2(5)$ and $V\downarrow R_i = U^l$, a direct sum of $l$ copies of an irreducible $KR_i$-submodule $U$ of dimension $2$.
\end{lem}

\pf  Suppose $R:=R_i \triangleleft G$. By Lemma \ref{norm1}, $V\downarrow R = U^l$ with $U$ irreducible and $(G_U)^U$ $\frac{1}{2}$-transitive. If $(R,\dim U)= (SL_2(5),2)$ then the conclusion holds, so suppose this is not the case. 
If $R^U$ is semiregular then $R$ is a Frobenius complement, so $R \cong SL_2(5)$; but then $\dim V$ must be 2 by Proposition \ref{sl25}(i), which we have assumed not to be the case.
Therefore $R^U$ is not semiregular, and so $U=V$ by Lemma \ref{norm1}(iv). In particular $F^*(G_0) = ZR$.

At this point we wish to apply \cite[Theorem 2.2]{KP}: this states that, with specified exceptions, any $p'$-subgroup of $GL_d(p)$ that has a normal irreducible quasisimple subgroup, has a regular orbit on vectors. In order to apply this, we need to establish that our quasisimple normal subgroup $R$ of $G$ acts irreducibly on $V$, regarded as an $\F_pR$-module. To see this, we go back to the proof of Lemma \ref{norm1}, letting $N: = R \triangleleft G$. Taking $U'$ to be an irreducible $\F_pR$-submodule of $V$, that proof shows that $R$ is faithful on $U'$, and that $G_{U'}$ is $\frac{1}{2}$-transitive on $U'$. Hence by the minimality of $G$, either $U'=V$ (which is the conclusion we want), or $G_{U'}^{U'}$ is semiregular or as in (ii) or (iii) of Theorem \ref{hope}. In the semiregular case, Proposition \ref{frob} implies that $R = SL_2(5)$ and $U'$ is a 2-dimensional $R$-module over some extension $K$ of $\F_p$, and this holds in (ii) and (iii) of Theorem \ref{hope} as well. However this can only happen if $\dim_KV = 2$ , contradicting our assumption that $(R,\dim U)\ne (SL_2(5),2)$. Hence $U'=V$, as desired.

Now we apply \cite[Theorem 2.2]{KP} which  determines all the possibilities for $G$ not having a regular orbit on $V$; these are 
\begin{itemize}
\item[(1)] the case with $R = A_c$ ($c<p$) and $V$ the deleted permutation module of dimension $c-1$, and
\item[(2)] the cases listed in Table \ref{kpegs}.
\end{itemize}

\vspace{2mm}
\no{\bf  Case (1) } In this case $G = Z_0H$ where $Z_0$ is a group of scalars and $H = A_c$ or $S_c$, and $V = \{(\a_1,\ldots ,\a_c) \in \F_p^c : \sum \a_i = 0\}$. If $v_1 = (1,-1,0,\ldots ,0)$ and $v_2 = (1,1,-2,0,\ldots ,0)$, one checks that the sizes of the $G$-orbits containing $v_1$ and $v_2$ are $\frac{c(c-1)|Z_0|}{(2,|Z_0|)}$ and $3|Z_0|{c \choose 3}$ respectively. These are not equal for any $c\ge 5$, contradicting $\frac{1}{2}$-transitivity.

\vspace{2mm}
\no {\bf Case (2) } In the case where $G/Z = U_4(2)$ and $(n,q) = (4,7)$, $G$ has two orbits on 1-spaces of sizes 40 and 360 (see \cite{liebaff}), and so cannot be $\frac{1}{2}$-transitive on $V^\sharp$. 
In each other case in Table 1, \cite[Theorem 2.2]{KP} gives the existence of a vector $v$ with stabiliser $G_v$ contained in a subgroup as indicated in column 4 of the table; and examination of the corresponding Brauer character of $G$ of degree $n$ in 
\cite{Atlas} gives the existence of another vector $u$ with stabiliser $G_u$ containing an element of order $m$, 
as indicated in column 5.  It follows in all cases that $G$ is not $\frac{1}{2}$-transitive. \hal

\begin{table}\label{kpegs}
\caption{Groups in case (2) of the proof of Lemma \ref{rinorm}}
\[
\begin{array}{|l|l|l|l|l|}
\hline
G/Z & n & q & G_v \le  & m \\
\hline
A_5 & 3 & 11 &  C_2 & 3 \\
S_5 & 4 & 7 & C_2 & 3 \\
S_6 & 5 & 7  & C_2 & 5 \\
A_6.2 & 4 & 7 & C_3 & 2 \\
A_6 & 3 & 19,31 & C_2,C_2 & 5,3 \\
A_7 & 4 & 11 & C_3 & 7 \\
L_2(7) & 3 & 11 & C_2 & 3 \\
L_2(7).2 & 3 & 25 & C_2 & 3 \\
U_3(3).2 & 7 & 5 & C_2 & 7 \\
U_3(3).2 & 6 & 5 & S_3 & 4 \\
U_4(2) & 4 & 7 & -- & -- \\
U_4(2) & 5 & 7,13,19 & S_4,V_4,C_2 & 5,5,5 \\
U_4(2).2 & 6 & 7,11,13 & D_{12},V_4,C_2 & 5,5,5 \\
U_4(2) & 4 & 13,19,31,37 & [18],[9],C_3,C_2 & 4,2,2,3 \\
 U_4(3).2 & 6 & 13,19,31,37 & W(B_3),S_3\times C_2,V_4,C_2 & 5,5,5,5 \\
U_5(2) & 10 & 7 & V_4 & 3 \\
Sp_6(2) & 7 & 11,13,17,19 & C_2^3,V_4,C_2,C_2 & 7,7,7,7 \\
\O_8^+(2) & 8 & 11,13,17,19,23 & W(B_3),S_4,S_3,V_4,C_2 & 7,7,7,7,7 \\
J_2 & 6 & 11 & S_3 & 4 \\
\hline
\end{array}
\]
\end{table}

\begin{lem}\label{k1}
We have $k=1$.
\end{lem}

\pf Suppose $k>1$. Assume first that $R_i \triangleleft G$ for all $i$. Then $N:=R_1R_2 \triangleleft G$; moreover $N$ is not a Frobenius complement by Proposition \ref{frob}, so is not semiregular on $V^\sharp$, and hence Lemma \ref{norm1}(iv) shows that $N$ is irreducible on $V$. Now Lemma \ref{rinorm} implies 
that 
\[
N = R_1R_2 = SL_2(5) \otimes SL_2(5) \le G \le \G L_4(q).
\]
Let $V = U\otimes W$ be a tensor decomposition preserved by $N$, with $\dim U = \dim W = 2$. 
If $q \ne p$ or $p^2$ with $p\le 61$, and also $q \ne 7^4$, then Proposition \ref{sl25} shows that the group induced by $G/Z$ on 1-spaces in $U$ has a regular orbit, and the same for $W$. Pick $\la u\ra$ and $\la w \ra$ in such orbits ($u\in U, w\in W$). Then $G_{\la u \otimes w\ra} \le Z$ and so $G_{u\otimes w} = 1$. Hence $G$ has a regular orbit on $V^\sharp$, a contradiction. 
And if $q = p$, $p^2$ or $7^4$, then  
\[
G \le  Z\cdot (SL_2(5) \otimes SL_2(5)).a  = Z\cdot R_1R_2.a \le \G  L_4(q),
\]
where $a$ divides 4. Here $G_0 = Z\cdot R_1R_2$. 
 Let $u_1,u_2$ be a basis of $U$ and $w_1,w_2$ a basis of $W$. Writing matrices relative to these bases, 
define $R_2^T = \{A^T: A \in R_2\}$. Then by \cite[Lemma 4.3]{arith}, for the vector $v = u_1\otimes w_1+u_2\otimes w_2$ we have
\begin{equation}\label{k2}
(G_0)_v = \{B \otimes B^{-T} : B \in R_1 \cap R_2^T\}.
\end{equation}
There is only one conjugacy class of subgroups $SL_2(5)$ in $GL_2(q)$, so we can choose bases $u_i,w_i$ such that $R_1 = R_2^T$; then for the corresponding vector $v$ the order of  $(G_0)_v$ is divisible by 60. On the other hand there are bases for which $R_1\cap R_2^T$ has order dividing 20, giving a vector stabilizer in $G$ of order coprime to 3. This contradicts $\frac{1}{2}$-transitivity.  

Thus not all the $R_i$ are normal subgroups of $G$. Relabelling, we may therefore take it that $G$ permutes $l$ factors $R_1,\ldots ,R_l$ transitively by conjugation, where $l>1$. Let $N = R_1\ldots R_l$. Lemma \ref{norm1}(iv) implies that $N$ is irreducible on $V$, so that $k=l$ and $F^*(G_0) = ZN$.  Now \cite[(3.16), (3.17)]{asch} implies that $N$ preserves a tensor decomposition $V = V_1\otimes \cdots \otimes V_k$ with $\dim V_i$ independent of $i$, $N \le \bigotimes GL(V_i)$ and $G \le N_{\G L(V)}(\bigotimes GL(V_i)) = (GL(V_1) \circ \cdots \circ GL(V_k)).S_k.\la \s \ra$ with $\s$ a field automorphism acting on all factors.

Let $G_1$ be the kernel of the natural map from $G$ to $S_k$, so that $G_1 = G \cap B$ where $B = 
 (GL(V_1) \circ \cdots \circ GL(V_k)).\la \s \ra$. There is a map $\phi: G_1 \rightarrow P\G L(V_1)$ which has image normalizing the simple irreducible group $T:=R_1/Z(R_1)$. 

Just as in the second paragraph of the proof of Lemma \ref{rinorm}, $N$ acts irreducibly on $V$, regarded as an $\F_pN$-module. 
It follows that $R_1$ acts irreducibly on $V_1$, regarded as an $\F_pR_1$-module: for if $W_1$ is a proper nonzero $\F_pR_1$-submodule of $V_1$, then by the transitivity of $G$ on the $R_i$, there is a proper nonzero $\F_pR_i$ submodule $W_i$ of $V_i$ for each $i$, and then $W_1\otimes \cdots \otimes W_l$ is an $\F_pN$-submodule of $V$, contradicting the $\F_pN$-irreducibility of $V$.

As in the proof of Lemma \ref{rinorm}, this means that we can apply \cite[Theorem 2.2]{KP} to the action of $G_1\phi$ on $V_1$. This shows that one of the following holds:
\begin{itemize}
\item[(a)]  $G_1\phi$ has a regular orbit on the 1-spaces of $V_1$;
\item[(b)] $T$ and $V_1$ are among the exceptions indicated in (1) and (2) of the proof of Lemma \ref{rinorm};
\item[(c)] $(T,\dim V_1) = (A_5,2)$.
\end{itemize}
\no Assume first that (a) holds and (c) does not. So $G_1\phi$ has a regular orbit on 1-spaces in $V_1$. Let $\la v \ra$ be a 1-space in such an orbit. Write also $v$ for the corresponding vector in the other $V_i$, and let $H$ be the stabiliser $(G_1)_{v\otimes \cdots \otimes v}$. Then $H$ fixes the 1-space $\la v \ra \otimes \cdots \otimes \la v \ra$, so by the choice of $v$, we have $H \le Z$, the group of scalars in $G$. Hence in fact $H=1$. It follows that $G_{v\otimes \cdots \otimes v}$ has order dividing $k!$.
Also, assuming $R_i \not \cong SL_2(r)$, there is an involution $r_i \in R_i\backslash Z$ fixing a nonzero vector $u_i \in V_i$, and hence we see that $G_{u_1\otimes \cdots \otimes u_k}$ has order divisible by $2^k$. 
However $2^k$ does not divide $k!$ so this is impossible. For $R_i \cong SL_2(r)$ we have $\dim V_i > 2$ (as we are assuming (c)  does not hold), and use a similar argument with an element of order 3 fixing a vector (which can be seen to exist from the character table of $SL_2(r)$ in \cite{dorn}).

Now consider case (b), where $T,V_1$ are as in (1) or (2) of the proof of Lemma \ref{rinorm}. For $T,V_1$ as in Table \ref{kpegs} (apart from $U_4(2)$ in dimension 4), let $v,u \in V_1$ be as in the last paragraph of the proof of Lemma \ref{rinorm}, and let $C$ be the group in the fourth column of Table \ref{kpegs} and $m$ the integer in the fifth. Then $(G_1)_{v\otimes \cdots \otimes v}$ is isomorphic to a subgroup of $C^k$, so that  $G_{v\otimes \cdots \otimes v}$ has order dividing $|C|^kk!$. 
On the other hand $(G_1)_{u\otimes \cdots \otimes u}$ has order divisible by $m^k$.  Since $G$ is $\frac{1}{2}$-transitive, this implies that $m^k$ divides $|C|^kk!$, which is not the case. 

The remaining cases in (b) are: $T=A_c\,(c<p)$, $V_1$ the deleted permutation module; and $T = U_4(2)$, $V_1 = V_4(7)$. 
In the latter case $T$ has two orbits on 1-spaces in $V_1$ with stabilizers of orders 72 and 648; so as above $G$ has a vector stabiliser of order dividing $72^kk!$ and another of order divisible by $648^{k-1}$, a contradiction. Now suppose $T=A_c\,(c<p)$ and $V_1$ is the deleted permutation module, which we represent as $\{(x_1,\ldots ,x_c) \in \F_p^c : \sum x_i=0\}$. By Bertrand's Postulate (see \cite{HW}) we can choose a prime $r$ such that $\frac{c}{2} < r < c$. Let $v_1,v_2$ be the following vectors in $V_1$:
\[
v_1 = (1^r,-r,0^{c-r-1}),\;\;v_2 = (1^{r-1},1-r,0^{c-r}).
\]
Then $G_{v_1\otimes \cdots \otimes v_1}$ has order divisible by $r^k$, while $G_{v_2\otimes \cdots \otimes v_2}$ has order dividing $m^kk!$, where $m = (r-1)!(c-r)!$ (note that $1-r\ne 1$ in $\F_p$, since $p>c$). Hence $r^k$ divides $k!$, a contradiction.

Finally consider case (c). Here $\dim V_i = 2$ and $R_i \cong SL_2(5)$; this case requires a special argument. Since $R_1$ is $\F_p$-irreducible on $V_1$, we must have $q=p$ or $p^2$, and hence $G \le Z\cdot (SL_2(5) \otimes \cdots \otimes SL_2(5)). S_k.\la \s\ra$ with $\s$ of order 1 or 2. Write $s = [\frac{k}{2}]$. As in the argument after (\ref{k2}), there is a vector $v \in V_1\otimes V_2$ whose stabilizer in $SL_2(5)\otimes SL_2(5)$ contains a diagonal copy of $SL_2(5)$. Tensoring $v$ with the corresponding vectors in $V_3\otimes V_4,\ldots , V_{2s-1}\otimes V_{2s}$ (and a further vector in $V_k$ if $k$ is odd), we see that there is a vector in $V$ with stabilizer in $G$ of order divisible by $60^s$. On the other hand there is a 1-space $\la w\ra$ in $V_1$ with stabilizer in $SL_2(5)/Z(SL_2(5))$ of order dividing 2, 3 or 5. Then $|G_{w\otimes \cdots \otimes w}|$ divides $t^kk!|\s|$ for some $t \in \{2,3,5\}$. Thus $60^{[k/2]}$ divides $t^kk!|\s|$. This is impossible unless $k$ is odd, $t=5$ and there is no 1-space in $V_1$ with stabilizer of order dividing 2 or 3. The latter can only hold if $q \equiv 3\hbox{ mod }4$ and  $q \equiv 2\hbox{ mod }3$. This implies that $q=p$ and $\s=1$, so that $60^{(k-1)/2}$ divides $5^kk!$. In particular $2^{k-1}$ divides $k!$, which is a contradiction for $k$ odd. \hal

\vspace{4mm}
We can now complete the proof of Theorem \ref{hope}. By Lemmas \ref{rinorm} and \ref{k1}, we have $F^*(G_0) = ZER_1$ where $R_1 = SL_2(5)$ and $E = O_2(G_0)$. Note that $p>5$ since $G$ is a $p'$-group, and so Lemma \ref{three} shows that 
$|E/\Phi(E)| \le 2^2$. Also by Lemma \ref{rinorm} we have $V\downarrow R_i = U^l$, a direct sum of $l$ copies of an irreducible $KR_i$-submodule $U$ of dimension $2$.

Suppose $E \not \le Z$, so that $E/\Phi(E) = 2^2$. Write $N = F^*(G_0)$. Proposition \ref{frob} shows that $N$ is not a Frobenius complement; hence 
Lemma \ref{norm1} shows that  $N$ is irreducible on $V$. Let $W$ be an irreducible $KE$-submodule of $V$. By Lemma \ref{norm1}, $E$ is faithful on $W$ (so $\dim W = 2$)  and  $G_W^W$ is a soluble $\frac{1}{2}$-transitive group. Such groups are classified in \cite[Theorem B]{pass2}. From this it follows that one of the following holds:
\begin{itemize}
\item[(a)] $G_W^W$ is a Frobenius complement (so $E$ is generalised quaternion);
\item[(b)] relative to some basis of $W$ we have $G_W^W = S_0(q)$, the group of monomial $2\times 2$ matrices of determinant $\pm 1$;
\item[(c)] $|W| = p^2$ with $p\in \{7,11,17\}$.
\end{itemize}
In case (c), $q=p$; also $p\ne 7,\,17$ as $SL_2(5) \not \le GL_2(p)$ for these values. Hence $V = U\otimes W = V_4(p)$ 
with $p = 11$, and a Magma computation shows that there is no such $\frac{1}{2}$-transitive group $G$ in this case.  
 
In  case (a), $G_W^W \le Z\cdot SL_2(3) < GL_2(q)$; and in (b), $G_W^W = Z\cdot 2^2 < Z\cdot SL_2(3).2 < GL_2(q)$.  
In either case it follows that $V = U \otimes W$ and $G \le Z\cdot (SL_2(5) \otimes (SL_2(3).2)) < GL_2(q) \otimes GL_2(q)< GL_4(q)$. Write $\bar G = GZ/Z$, so that $\bar G \le A_5\times S_4$.

We saw in the proof of Proposition \ref{sl25} that at least $q-62$ of the elements of $P_1(U)$ lie in regular orbits of $A_5$. 
Similarly, at least $q-32$ elements of $P_1(W)$ lie in regular orbits of $S_4$. Hence if $q>61$ then, picking $\la u\ra \in P_1(U)$ and $\la w \ra \in P_1(W)$ in regular orbits, we see that $u\otimes w$ lies in a regular orbit of $G$ on $V^\sharp$. This is a contradiction, since $G$ is $\frac{1}{2}$-transitive but not semiregular. 
Hence $q\le 61$. Now a Magma computation shows that no $\frac{1}{2}$-transitive groups arise in cases (a) and (b) as well. 

Thus we finally have $F^*(G_0) = ZR_1$ with $R_1 = SL_2(5)$ and $V\downarrow R_1 = U^l$, $\dim U = 2$.
Here $G/Z$ is $A_5$ or $S_5$, so $l=1$. Now Proposition \ref{sl25}(iii) shows that $q=11$, 19,  29 or 169 and $G$ is as in conclusion (ii) or (iii) of Theorem \ref{hope}. This is our final contradiction to the assumption that $G$ is a minimal counterexample.

\vspace{2mm}
This completes the proof of Theorem \ref{hope}. 

\section{Proof of Proposition \ref{khalf}}

Let $k\ge 2$ and suppose that $X$ is a $(k+\frac{1}{2})$-transitive permutation group of degree $n$. Assume that $X$ is not $k$-trransitive. We refer to \cite[\S 2]{K} for the list of 2-transitive groups, and to \cite[\S 7.6]{DM} for a discussion of sharply $k$-transitive groups.

The proposition is trivial if $X$ is $A_n$ or $S_n$, so assume this is not the case. Then $k\le 5$, as there are no 6-transitive groups apart from $A_n$ and $S_n$. Apart from $A_n$ and $S_n$, the only 5-transitive groups are the Mathieu groups $M_{12}$ and $M_{24}$, and the only 4-transitive, not 5-transitive, groups are $M_{11}$ and $M_{23}$. The groups $M_{11}$ and $M_{12}$ are sharply 4- and 5-transitive respectively; and in $M_{23}$, a 4-point stabilizer has orbits of size 3 and 16, so that $M_{23}$ is not $4\frac{1}{2}$-transitive and also $M_{24}$ is not $5\frac{1}{2}$-transitive. This gives the proposition for $k=4$ or 5.

Next let $k=3$. Then $X$ is a 3-transitive but not 4-transitive group, hence is one of the following: $AGL_d(2)$ (degree $2^d$); 
$2^4.A_7$ (degree $2^4$); $M_{11}$ (degree 12); $M_{22}$ or $M_{22}.2$ (degree 22); or a 3-transitive subgroup of $P\Gamma L_2(q)$ (degree $q+1$). The affine groups here are not $3\frac{1}{2}$-transitive, as a 3-point stabilizer fixes a further point. Neither are $M_{11}$, $M_{22}$ or $M_{22}.2$ as 3-point stabilizers have orbits of size 3,6 or 3,16. Finally, 
suppose that $X$ is a 3-transitive subgroup of $P\Gamma L_2(q)$. There are two possible sharply 3-transitive groups here, namely $PGL_2(q)$ and a group $M(q_0^2):=L_2(q_0^2).2$ with $q=q_0^2$ and $q$ odd, which is an extension of $L_2(q_0^2)$ by a product of a diagonal and a field automorphism. Assuming that $X$ is not one of these, it must be the case that a 3-point stabilizer $X_{\a\b\g} = \la \phi \ra$, where $\phi$ is a field automorphism. Since $X$ is  $3\frac{1}{2}$-transitive, $\la \phi \ra$ acts semiregularly on the remaining $q-2$ points, so any nontrivial power of $\phi$ must fix exactly 3 points. It follows that $q=2^p$ with $p$ prime,  and  $\phi$ has order $p$, which is the example in conclusion (iii) of Proposition \ref{khalf}.

Now suppose that $k=2$. Consider first the case where $X$ is almost simple, and let $T = {\rm soc}(X)$. When $T$ is not $L_2(q)$, $Sz(q)$ or $^2\!G_2(q)$, the arguments in \cite[\S 3]{K} show that a 2-point stabilizer $X_{\a\b}$ has orbits of unequal sizes on the remaining points, contradicting $2\frac{1}{2}$-transitivity. The groups with socle $L_2(q)$ are in conclusion (iv) of Proposition \ref{khalf}. 
If $T=\,^2\!G_2(q)$ (of degree $q^3+1$), then $X_{\a\b}$ has order $(q-1)f$, where $f=|X:T|$ is odd, and $X_{\a\b}$  is generated by an element $x$ of order $q-1$ and a field automorphism of odd order $f$. This group has a unique involution $x^{(q-1)/2}$ which fixes $q+1$ points. It follows that some nontrivial orbits of $X_{\a\b}$ have odd size and some have even size, contrary to 
$2\frac{1}{2}$-transitivity. Now consider $T = Sz(q)$, of degree $q^2+1$. If $X=T$ then it is a Zassenhaus group, and is in (iv) of the proposition. Otherwise, $X = \la T,\phi\ra$ where $\phi$ is a field automorphism of odd order $f$, say, and $\phi$ fixes $q_0^2+1$ points, where $q = q_0^f$. For suitable $\a,\b$ we have $X_{\a\b} = \la x,\phi\ra$, where $x$ has order $q-1$, and $\la x\ra$ has $q+1$  orbits of size $q-1$. Now $\phi$ fixes points in some of these orbits, so by $2\frac{1}{2}$-transitivity it must fix a point in each of them. But $|{\rm fix}(\phi| = q_0^2+1 < q+1$, which is a contradiction. 

Finally, suppose $X$ is affine (with $k=2$). Write $X = T(V)X_0 \le AGL(V)$, where $n=|V|$, $T(V)$ is the translation subgroup, and $X_0 \le GL(V)$. We refer to \cite[\S 2(B)]{K} for the list of possibliities for the transitive linear group $X_0$. If $X_0 \triangleright SL_d(q)\,(n=q^d, d\ge 2)$, $Sp_d(q)'\,(n=q^d, d\ge 4)$ or $G_2(q)'\,(n=q^6)$, the arguments in \cite[\S 4]{K} show that for some $v \in V^\sharp$, $X_{0v}$ has nontrivial orbits of unequal sizes. In cases (6-8) of  \cite[\S 2(B)]{K}, we have $X_0 \triangleright SL_2(5)$, $SL_2(3)$, $2^{1+4}$ or $SL_2(13)$, and $n \in \{3^4,3^6,5^2,7^2,11^2,19^2,23^2,29^2,59^2\}$; in each case $n-2$ is coprime to the order of a 2-point stabilizer $X_{0v}$, so it follows by  $2\frac{1}{2}$-transitivity that $X_{0v} = 1$. In other words, $X$ must be sharply 2-transitive, as in conclusion (ii) of the proposition.

It remains to deal with the case where $X \le A:= A\Gamma L_1(q)$ ($n=q$). Here $A_{01}$ consists of field automorphisms, so if we pick $v \in \F_q$ such that $v$ lies in no proper subfield of $\F_q$, then $A_{01v} = 1$. Hence by  $2\frac{1}{2}$-transitivity, all 3-point stabilizers in $X$ are trivial -- that is, $X$ is a Zassenhaus group. It is well known that the non-sharply 2-transitive Zassenhaus groups in the 1-dimensional affine case are just $A\G L_1(2^p)$ with $p$ prime, as in (iv) of the proposition. This is easy to see: we have $X_{01} = \la \phi\ra$, where $\phi$ is a field automorphism, and this acts semiregularly on $\F_q\setminus \{0,1\}$; hence, as argued at the end of the $k=3$ case above, $q=2^p$ with $p$ prime and $X = A\G L_1(2^p)$, as required.

This completes the proof of Proposition \ref{khalf}.

\vspace{1cm}

\no Martin W. Liebeck, Dept. of Mathematics, Imperial College, London SW7 2BZ, UK, 
email: m.liebeck@imperial.ac.uk

\vspace{4mm}
\no Cheryl E. Praeger, School of Mathematics and Statistics, University of Western Australia, Western Australia 6009, 
email: praeger@maths.uwa.edu.au

\vspace{4mm}
\no Jan Saxl, DPMMS, CMS, University of Cambridge, Wilberforce Road, Cambridge CB3 0WB, UK, 
email: saxl@dpmms.cam.ac.uk
\end{document}